\documentclass[12pt, leqno,a4paper]{article}
\usepackage{latexsym}
\usepackage{amscd}
\usepackage{amssymb}
\def\noproof{{\unskip\nobreak\hfill\penalty50\hskip2em\hbox{}%
     \nobreak\hfill$\Box$\parfillskip=0pt%
     \finalhyphendemerits=0\par}}
\def\enddemo{\ifmmode\eqno\Box\else\noproof\vskip0.8truecm\fi}
\newtheorem{theorem}{Theorem}[section]

\newtheorem{proposition}[theorem]{Proposition}
\newtheorem{corollary}[theorem]{Corollary}

\newtheorem{remark}[theorem]{Remark}
\newtheorem{lemma}[theorem]{Lemma}
\newcommand{\Spec}{{\mbox{Spec}\,}}
\newcommand{\lra}{\longrightarrow}
\newcommand{\Hom}{\mbox{Hom}}
\newcommand{\Ext}{\mbox{Ext}}

\newcommand{\Tor}{{\mbox{Tor}}}
\newcommand{\Pic}{{\mbox{Pic}}}
\newcommand{\ltor}{{\ell\mbox{\scriptsize{-tor}}}}

\newcommand{\pt}{{\pi_1^{t,ab}}}

\newcommand{\cor}{{\mbox{cor}}}
\newcommand{\fX}{{\frak{X}}}

\newcommand{\cZ}{{\cal Z}}

\newcommand{\cF}{{\cal F}}

\newcommand{\cR}{{\cal R}}
\newcommand{\bZ}{{\Bbb Z}}
\newcommand{\bZt}{{\bZ_{tr}}}
\newcommand{\bQ}{{\Bbb Q}}
\newcommand{\bA}{{\Bbb A}}
\newcommand{\bF}{{\Bbb F}}
\newcommand{\bL}{{\Bbb L}}
\newcommand{\io}{{\iota}}
\newcommand{\et}{{\mbox{\scriptsize{\'et}}}}
\newcommand{\se}{{\subseteq}}

\begin{document}

\title{On the Albanese map for smooth quasi-projective varieties}
\author{By Michael Spie{\ss} and Tam\'as Szamuely}
\date{}
\maketitle

\section{Introduction}

Consider an algebraically closed field $k$ of characteristic $p\geq 0$ and a smooth quasi-projective $k$-variety $X$. If $X$ is in fact projective, a famous theorem due to A. A. Roitman (\cite{roitman}, see also \cite{bloch}) asserts that the Albanese map
\begin{equation}
\label{eqn:alb}
alb_X: {CH_0(X)}^0\lra Alb_X(k)
\end{equation}
from the Chow group of zero-cycles of degree 0 on $X$ to the group of $k$-points of the Albanese variety induces an isomorphism on prime-to-$p$ torsion subgroups (later J. S. Milne proved that the isomorphism holds for $p$-primary torsion subgroups as well, cf. \cite{milne}). As a well-known counter-example of Mumford shows, one cannot expect the map $alb_X$ itself to be an isomorphism in general. Still, Kato and Saito (\cite{kasa}, Section 10) have established the bijectivity of $alb_X$ in the case when $k$ is the algebraic closure of a finite field (in fact, in this case both groups are torsion). Moreover, bijectivity over $k={\overline \bQ}$ has been conjectured by Bloch and Beilinson, as a consequence of some expected standard features of the conjectural category of mixed motives over ${\overline \bQ}$.

In the present paper we generalise the theorems of Roitman and Kato/Sai\-to to the case when $X$ is not necessarily projective but admits a smooth compactification. Here the natural target for the Albanese map is the generalised Albanese variety introduced by Serre \cite{serre}. If $X$ is a curve, this variety is a generalised Jacobian in the sense of Rosenlicht \cite{ros} and for $X$ proper it is the usual Albanese. In general, it is a semi-abelian variety universal for morphisms of $X$ into semi-abelian varieties; it is related to the Picard variety by a duality theorem (see sections 3 and 4 for more details). The generalisation of $alb_X$ to this context is a map 
\begin{equation}
\label{eqn:galb}
alb_X: {h_0(X)}^0 \lra Alb_X(k)
\end{equation}
where the group on the left is the degree zero part of Suslin's 0-th algebraic singular homology group defined in \cite{sv}; it coincides with $CH_0(X)^0$ for $X$ proper. The map (\ref{eqn:galb}) was first constructed by N. Ramachandran in \cite{ram}; we give a simpler proof for the ``reciprocity law'' implying its existence in Section 3. 
  
Our results are as follows.

\begin{theorem} 
\label{theorem:thm1}
Let $k$ be an algebraically closed field of characteristic $p\geq 0$ and let $X$ be a smooth quasi-projective variety over $k$. Assume that there exists a smooth projective connected $k$-variety $\fX$  containing $X$ as an open subscheme. Then the Albanese map (\ref{eqn:galb}) induces an isomorphism on prime-to-$p$ torsion subgroups.
\end{theorem}

Note that the required smooth compactification $\fX$ exists if $k$ is of characteristic 0 or if $X$ is of dimension $\leq 3$ and $p\geq 5$, by virtue of the desingularisation theorems of Hironaka and Abhyankar.

\begin{theorem} 
\label{theorem:thm2} We keep the hypotheses of Theorem \ref{theorem:thm1} and assume moreover that $k$ is the algebraic closure of a finite field. Then (\ref{eqn:galb}) is an isomorphism of torsion groups of finite corank.
\end{theorem}

Our method for proving Theorem \ref{theorem:thm1} is new even in the proper case and is (at least to our feeling) more conceptual than the previous ones. One of the key observations is that thanks to its functoriality and homotopy invariance properties, the (generalised) Albanese variety can be regarded as an object of Voevodsky's triangulated category of effective motivic complexes $DM^{eff}_-(k)$ and in fact for smooth varieties the Albanese map can be interpreted as a morphism in this category. As we shall see in Section 5, the theorem then results from the fundamental isomorphism that relates the algebraic singular cohomology group $h^1(X,\bZ/n\bZ)$ to the \'etale cohomology group $H^1_{\et}(X,\bZ/n\bZ)$ according to Suslin and Voevodsky \cite{sv}. The proof of Theorem \ref{theorem:thm2} is a direct generalisation of the argument given in (\cite{kasa}, section 10), using the ``tamely ramified class field theory'' developped in  \cite{ss}.

Finally it should be mentioned that during recent years fruitful efforts have been made for generalising Roitman's theorem to singular complex projective varieties (see \cite{bs} and the references quoted there). Our generalisation seems to be unrelated to this theory except perhaps in the case when $X$ is the complement of the singular locus of a complex projective variety.

A word on notation:  For an abelian group $A$ and a nonzero integer $n$ we denote by $A_n$ the $n$-torsion subgroup of $A$ and we write $A/n$ as a shorthand for $A/nA$. For a prime number $\ell$ we let $A_{\ltor}$ be the $\ell$-primary component of the torsion subgroup of $A$.
 
\section{Review of Algebraic Singular Homology}

In this section and the next, we discuss the definition of the map (\ref{eqn:galb}) and the groups involved. We shall work over an arbitrary perfect field $k$ and $X$ will be any $k$-variety (i.e. an integral separated $k$-scheme of finite type over $k$). 

Let us recall the definition of the algebraic singular homology groups introduced in \cite{sv}. Consider for an integer $n\geq 1$ the algebraic $n$-simplex $\Delta^n= \Spec k[T_0, \ldots, T_n]/(T_0 + \ldots + T_n - 1)$. Denote by $C_n(X)$ the free abelian group generated by those integral closed subschemes $Z$ of $X\times\Delta^n$ for which the projection $Z\to \Delta^n$ is finite and surjective. Any monotonous nondecreasing map $\alpha: \{0,1,\ldots, m\} \to \{0,1,\ldots, n\}$ induces a morphism $\Delta^m \to \Delta^n$ and thus a homomorphism $\alpha^*: C_n(X) \to C_m(X)$ via pull-back. These maps endow the set of the $C_{n}(X)$ with the structure of a simplicial abelian group; we denote by $C_{\bullet}(X)$ the associated chain complex. For an abelian group $A$ the $n$-th algebraic singular homology group $h_n(X, A)$ of $X$ with coefficients in $A$ is defined as the $n$-th homology of the complex $C_{\bullet}(X)\otimes A$ and the $n$-th algebraic singular cohomology $h^n(X, A)$ as the $n$-th cohomology of $\Hom(C_{\bullet}(X, \bZ), A)$.  For $A=\bZ$ we shall simply write $h_n(X)$ for $h_n(X,\bZ)$ etc.

The group $h_0(X)$ has the following concrete description. Let $\cZ(X)$ be the free abelian group with basis the set $X_0$ of closed points of $X$. Then $h_0(X)$ is the quotient of $\cZ(X)$ by the submodule $\cR$ generated by $i_0^*(Z)-i_1^*(Z)$, where $i_{\nu}: X \to X\times \bA^1, x\mapsto (x, \nu)$, $\nu=0,1$ and $Z$ runs through all closed integral subschemes of $X\times \bA^1$ such that the projection $Z \to \bA^1$ is finite and surjective. There is a natural degree map $\cZ(X) \lra \bZ$ given by the formula
\[
\sum_i n_i P_i \mapsto \sum_i n_i [k(P_i):k],
\]
of which we denote the kernel by $\cZ(X)^0$. Using the fact that the projections $Z \to \bA^1$ are finite and flat, one checks that $Z(X)^0$ contains $\cR$; the quotient $\cZ(X)^0/ \cR$ will be denoted by $h_0(X)^0$.

For the proofs we shall also need the sheafified version of the above construction. For this, denote by $Sm/k$ the category of smooth schemes of finite type over $k$. Let $\cF$ be an abelian presheaf on $Sm/k$, i.e. a contravariant functor from $Sm/k$ to the category of abelian groups. For any $m\geq 0$ we may define a presheaf $\cF_m$ by the rule $\cF_m(X)={\cF}(X\times \Delta^m)$. Together with the operations induced from the cosimplicial scheme $\Delta^{\bullet}$ these presheaves assemble to form a simplicial presheaf whose associated chain complex we denote by $C_{\bullet}({\cal F})$. If $\cF$ is {\em homotopy invariant}, i.e. if $\cF(X\times \bA^1) = \cF(X)$ for all $X\in Sm/k$, then the natural augmentation map $C_{\bullet}(\cF) \to \cF$ given by the identity in degree 0 is a map of complexes and in fact a quasi-isomorphism (here we view $\cF$ as a complex concentrated in degree 0). Indeed, in view of the canonical isomorphism $\Delta^n\cong \bA^n$ in this case $C_{\bullet}(\cF)$ is none but the complex associated to the constant simplicial presheaf defined by $\cF$. 

We also recall the notion of {\em presheaves with transfers} from Section 2 of \cite{voe}. These are contravariant additive functors with values in abelian groups from  the category $SmCor(k)$ whose objects are smooth schemes of finite type over $k$ and where a morphism from an object $X$ to an object $Y$ is a {\em finite correspondence,} i.e. an element of the free abelian group $Cor(X,Y)$ generated by those integral closed subschemes $Z$ of $X\times Y$ for which the projection $Z\to X$ is finite and surjective over a component of $X$. ({\em Note:} In \cite{sv}, of which we shall use the results, presheaves with transfers are defined in a different way. However, the proof of the so-called rigidity theorem goes through with the above definition and this is sufficient for our applications.)  
Now the link with the algebraic singular complex is the following. For a separated $k$-scheme $X$ the rule $U \mapsto Cor(U, X)$ defines a presheaf with transfers on which we denote by $\bZt(X)$; actually it is a sheaf for the \'etale topology on $SmCor(k)$. Then by definition $C_{\bullet}(\bZt(X))(k) = C_{\bullet}(X)$. 
 
\section{The Generalised Albanese Map}

Keeping the assumptions of the previous section, we begin by recalling some basic facts about the generalised Albanese variety $Alb_X$ of $X$ which was introduced in \cite{serre}. Although it can be defined for an arbitrary variety (see \cite{ram}) we assume from now on for simplicity that $X$ has a $k$-rational point. Then $Alb_X$ is a semiabelian variety satisfying the following universal property: for every $k$-rational point $P$ of $X$ there is a morphism $\io_P:X \to Alb_X$ such that $\io_P(P) = 0$ and if $(B, f)$ is a pair consisting of a semiabelian variety $B$ and a morphism $f:X \to B$ mapping $P$ to $0_B$ there is a unique morphism $g: Alb_X \to B$ of group schemes with $g\circ \io_P = f$. If $X$ is proper, then $Alb_X$ is the Albanese variety in the classical sense. If $X$ is a curve, it coincides with Rosenlicht's generalised Jacobian for the modulus defined by the sum of points at infinity. 

The assignment $X \mapsto Alb_X$ is a covariant functor for arbitrary morphisms of varieties. 
Moreover, there is also a contravariant functoriality of $Alb_X$ with respect to finite flat morphisms $f:\,X\to Y$ which we now briefly explain. Assume first $k$ is algebraically closed. Then mapping a closed point $Q$ of $Y$ to the pull-back zero-cycle $f^*(Q)$ defines a morphism of $Y$ into the $d$-fold symmetric product $Sym^d(X)$, where $d$ is the degree of $f$. On the other hand, for a fixed $k$-point $P$ of $Y$ the zero-cycle $f^*(P)=P_1+\dots+P_d$ defines a morphism $Sym^d(X)\to Alb_X$ via the sum of the maps $\iota_{P_i}$ $(1\leq i\leq d)$. The composite of these two maps sends $P$ to 0 in $Alb_X$, hence by definition of $Alb_Y$ factors as the composite of $\iota_P$ with a morphism $f^*:Alb_Y\to Alb_X$ which is easily seen to be independent of $P$; it is the one we were looking for. For arbitrary perfect $k$ one remarks that carrying out the above construction over the algebraic closure gives a Galois-equivariant map, which is thus defined over $k$. 

We denote by 
\[
a_X: \cZ(X)^0 \lra Alb_X(k)
\]
the homomorphism 
\[
\sum_i n_i P_i \mapsto \sum_i n_i \cor_{k(P_i)/k}(\io_P(P_i))
\]
(note that the sum is independent of the base point $P\in X(k)$ since we have $\io_{P_1}(P_3) = \io_{P_2}(P_3) + \io_{P_1}(P_2)$ for any two $k$-rational points $P_1, P_2$ of $X$ and closed point $P_3$). For a morphism $f:X \to Y$ of varieties the diagram
\begin{equation}
\label{eqn:cofun}
\begin{CD}
Z(X)^0 @> a_X >> Alb_X(k)\\
@VV f_* V @VV f_* V\\
Z(Y)^0 @> a_Y >> Alb_Y(k)\\
\end{CD}
\end{equation}
commutes where the vertical maps are induced by $f$ through covariant functoriality. Similarly, for a finite flat $f:\, Y\to X$ the diagram

\begin{equation}
\label{eqn:contra}
\begin{CD}
Z(Y)^0 @> a_Y >> Alb_Y(k)\\
@VV f^* V @VV f^* V\\
Z(X)^0 @> a_X >> Alb_X(k)\\
\end{CD}
\end{equation}
commutes.

Using these functoriality properties we can give an easy proof of the following ``reciprocity law'' which immediately yields the existence of the map $alb_X$ as in (\ref{eqn:galb}).

\begin{lemma}\label{recip}
With notations as above, the subgroup $\cR\subset\cZ(X)^0$ is contained in the kernel of the map $a_X$.
\end{lemma}

\noindent{\it Proof.} Let $Z\se X\times \bA^1$ be a closed integral subscheme such that the projection $q:Z \to \bA^1$ is finite and surjective (hence flat) and denote by $p:Z \to X$ the other projection. By the commutativity of (\ref{eqn:cofun}) and (\ref{eqn:contra}) we have 
\[
a_X(i_0^*(Z)-i_1^*(Z)) = a_X(p_*(q^*((0)-(1)))) = p_*(q^*(a_{\bA^1}((0)-(1)))).
\]
Since the generalised Albanese of $\bA^1$ is trivial (any map of $\bA^1$ into a semi-abelian variety being constant), it follows that the left hand side is $0$.
\enddemo

Just like the previous one, this section concludes with some remarks concerning the sheafification of the above construction. We only consider the case when $X$ is {\em smooth}. For our purposes it is convenient to replace $Alb_X$ by the ``Albanese scheme'' $\widetilde{Alb}_X$ introduced by Ramachandran in (\cite{ram}, section 2.2). It is a smooth commutative group scheme which is an extension of the constant group scheme $\bZ$ by $Alb_X$. The choice of a base point defines a splitting i.e. an isomorphism $\bZ\times Alb_X\cong\widetilde{Alb}_X$. Again $\widetilde{Alb}_X$ satisfies a universal property and we have a canonical morphism $\io: X \to \widetilde{Alb}_X$ which does not depend on the choice of a base point. The Albanese map can be naturally extended to a map $alb_X: h_0(X) \to  \widetilde{Alb}_X(k)$.

Consider now the abelian presheaf on $Sm/k$ represented by the group scheme $\widetilde{Alb}_X$ which we also denote by $\widetilde{Alb}_X$. It is a sheaf for the \'etale (even the $fppf$) topology. 

\begin{lemma}
\label{transfers}
The \'etale sheaf $\widetilde{Alb}_X$ is a homotopy invariant presheaf with transfers.
\end{lemma}

\noindent{\em Proof.}
Homotopy invariance is again a consequence of the fact that there is no non-constant map $\bA^1\to \widetilde{Alb}_X$. To construct transfer maps, we can work more generally with an arbitrary commutative group scheme $G$. Let $X,Y \in Sm/k$ and let $Z\subset X\times Y$ closed integral subscheme finite  and surjective over a component of $X$. As explained before Theorem 6.8 of \cite{sv}, to $X$ one can associate a canonical map $\alpha_Z:\,X\to Sym^d(Y)$, where $d$ is the degree of the projection $Z\to X$. Now given a map $Y\to G$, it induces a map $Sym^d(Y)\to Sym^d(G)$, whence we obtain the required map $X\to G$ by composing by $\alpha_Z$ on the left and by the summation map on the right.  
\enddemo

The lemma implies that there is a unique morphism of presheaves 
\begin{equation}\label{albtransfer}
\bZt(X) \to \widetilde{Alb}_X
\end{equation} 
which sends the correspondence associated to the identity map $id:\,X\to X$ to the Albanese map $\io\in\widetilde{Alb}_X(X)$. By applying the functor $C_{\bullet}(\,\,)$ we get a map  $C_{\bullet}(\bZt(X)) \to C_{\bullet}(\widetilde{Alb}_X)$. Composing it with the augmentation map on the right (existing by homotopy invariance of $\widetilde{Alb}_X$; see the previous section) yields the map of complexes of \'etale sheaves with transfers
\begin{equation}
\label{eqn:sgalb}
C_{\bullet}(\bZt(X)) \lra \widetilde{Alb}_X.
\end{equation}
Here we again consider $\widetilde{Alb}_X$ as a complex concentrated in degree 0. Since this is a morphism of complexes, it factors through the 0-th homology presheaf $H_0(C_{\bullet}(\bZt(X)))$; as $\widetilde{Alb}_X$ is an \'etale sheaf, it even factors through the associated \'etale sheaf $H_0(C_{\bullet}(\bZt(X)))_{\et}$. This means that the existence of the map (\ref{eqn:sgalb}) subsumes a sheafified version of the reciprocity law (perceptive readers have already noted the similarity of the argument with the proof of Lemma \ref{recip}). In fact, by passing to sections over $\Spec (k)$ in (\ref{eqn:sgalb}) and to homology we get the generalised Albanese map (\ref{eqn:galb}). Finally, we remark for later use that the map (\ref{albtransfer}) can be obtained as a composite of (\ref{eqn:sgalb}) with the natural morphism of complexes $\bZt(X)\to C_{\bullet}(\bZt(X))$ (again with $\bZt(X)$ placed in degree 0 on the left). 

\begin{remark}\rm
In the terminology of \cite{voe}, Lemma \ref{transfers} states that the sheaf $\widetilde{Alb}_X$ defines an object in the category $DM^{eff}_-(k)$ of effective motivic complexes; on the other hand, $C_{\bullet}(\bZt(X))$ is precisely the motivic complex that Voevodsky associates to the smooth variety $X$. Therefore the map (\ref{eqn:sgalb}), which was shown above to be a morphism in $DM^{eff}_-(k)$, can be regarded as the ``motivic interpretation'' of the Albanese map. 
\end{remark}

\section{Relation to Tame Abelian Covers}

Assume now we are in the situation of Theorem \ref{theorem:thm1}. In this case  $Alb_X$ has been described by Serre in his expos\'e \cite{serrebis} as an extension of the abelian variety $Alb_{\fX}$ by a torus whose rank is equal to the rank of the subgroup $B_X$ of divisors on $\fX$ which are algebraically equivalent to zero and whose support is contained in $\fX-X$. As a consequence of this result one gets that, just as in the proper case (see \cite{kala}, \cite{milne}), abelian \'etale covers of $X$ of degree prime to $p$ can be obtained by pull-back of those from $Alb_X$ ``up to a finite obstruction''. 

More precisely, denote by $S$ the finite torsion subgroup of the quotient $Q$ of the N\'eron-Severi group $NS(\fX)$ of $\fX$ by the subgroup generated by classes of divisors with support in $\fX-X$. Then we have:

\begin{proposition}
\label{proposition:prop1}
For every integer $n$ prime to $p$ there is an exact sequence\begin{equation}
\label{eqn:ses}
0 \lra \Hom({Alb_X(k)}_n, \bZ/n) \lra H^1_{\et}(X, \bZ/n) \lra \Tor(S, \bZ/n) \lra 0.
\end{equation}
\end{proposition}

\noindent{\em Proof.} Without loss of generality we may assume that either $X=\fX$ or $X$ is the complement of a divisor in $\fX$. Indeed, if not, then we may find an open subscheme $X'$ of $\fX$ containing $X$ which is the complement of a divisor in $\fX$ and such that the codimension of $X'-X$ in $X'$ is at least $2$. Then we have canonical isomorphisms $Alb_X\cong  Alb_{X'}$ (see \cite{ram}, Corollary 2.4.5) and $H^1_{\et}(X, \bZ/n)\cong H^1_{\et}(X', \bZ/n)$ (a consequence of Zariski-Nagata purity; see \cite{sga2}, expos\'e X) which are compatible with the Albanese map.

In the rest of the proof we use the language of 1-motives for which we refer the readers to Section 10 of \cite{deligne}. In Section 2.4 of \cite{ram}, Ramachandran defined a distinguished triangle in the derived category of bounded complexes of smooth commutative group schemes
\begin{equation}
\label{eqn:tri}
M^1(X) \lra M^{\tau}(X) \lra S[-1] \lra M^1(X)[1].
\end{equation}
Here $S[-1]$ is the complex $[0 \to S]$ concentrated in degrees 0 and 1, $S$ being considered as a finite constant group scheme. The term $M^1(X)$ is the ``Picard 1-motive'' of $X$ which is by definition the 1-motive whose term in degree 0 is the constant group scheme attached to the finitely generated free abelian group $B_X$ introduced in the beginning of this section and the degree 1 term is the Picard variety of $\fX$. Finally, to define the ``Picard $\tau$-motive'' $M^{\tau}(X)$, one observes first that if $I$ is the set of codimension 1 components of $\fX\setminus X$, the complex of group schemes $\bZ^I\to \Pic(\fX)]$ is naturally an extension of the complex $[\bZ^I/B_X \to NS(\fX)]\simeq [0 \to Q] = Q[-1]$ by the 1-motive $M^1(X)$; one then obtains $M^{\tau}(X)$ by pulling back this extension by the natural inclusion $S[-1]\to Q[-1]$. By passing to $k$-valued points in (\ref{eqn:tri}), tensoring with $\otimes^{\bL} \bZ/n$ and passing to cohomology we obtain the exact sequence
\[
0 \lra H^0(M^1(X)(k)\otimes^{\bL} \bZ/n) \lra H^0(M^{\tau}(X)(k)\otimes^{\bL} \bZ/n) \lra \Tor(S, \bZ/n) 
\]
\[
\lra H^1(M^1(X)(k)\otimes^{\bL} \bZ/n).
\]
We can identify this exact sequence to that of the proposition as follows. The last term vanishes by the $n$-divisibility of the group of rational points of a smooth connected commutative group scheme over our algebraically closed field $k$. According to (\cite{ram}, Theorem 2.4.6), the second term is canonically isomorphic to $H^1_{\et}(X, \mu_n)$, hence to $H^1_{\et}(X, \bZ/n)$ by choosing an isomorphism $\mu_n\cong\bZ/n$. Similarly, it will suffice to identify the first term with $\Hom({Alb_X(k)}_n, \mu_n)$. For this, recall that, as explained by (\cite{ram}, Theorem 2.4.4), the main result of \cite{serrebis} can be reinterpreted by saying that the Cartier dual of the 1-motive $M^1(X)$ is the 1-motive $[0 \to Alb_X]$; in particular, the toric part $T$ of $Alb_X$ has character group $B_X$. There is an exact sequence
\[
0 \lra \Hom(Alb_{\fX}(k)_n, \mu_n) \lra \Hom(Alb_X(k)_n, \mu_n) \lra \Hom(T_n, \mu_n) \lra 0
\]
coming from writing $Alb_X$ as a canonical extension of $Alb_{\fX}$ by $T$. On the other hand, the stupid filtration on the two-term complex $M^1(X)$ induces an exact sequence
\[
0 \lra\Tor(\Pic^0_\fX(k)_n, \bZ/n) \lra H^0(M_1(k)\otimes^{\bL} \bZ/n) \lra B_X/n \lra 0.
\]
In view of the fact that $\Tor(\Pic^0_\fX(k)_n, \bZ/n)\cong \Pic^0_\fX(k)_n$, which in turn is the dual group of $Alb_{\fX}(k)_n$ via the $e_m$-pairing, Serre's result implies that this sequence is canonically isomorphic to the previous one.
\enddemo

In the next section we shall use the proposition through the corollary:

\begin{corollary}
\label{corollary:cor1}
Under the assumptions of the proposition there is a canonical isomorphism
\[
\Hom(H^1_{\et}(X, \bZ_{\ell}),\bQ_{\ell}/\bZ_{\ell})\cong {Alb_X(k)}_{\ltor}
\]
for any prime number $\ell\neq p$.
\end{corollary}

\noindent{\em Proof.} Passing to the inverse limit by making $n$ run over powers of $\ell$ in (\ref{eqn:ses}) yields an exact sequence as the groups $\Tor(S, \bZ/\ell^n)$ are finite. In fact, their inverse limit is trivial, so the corollary follows by dualising.
\enddemo

In the remainder of this section, which will only be needed for the proof of Theorem \ref{theorem:thm2}, we strengthen the result of the proposition to obtain a description of the abelianised tame fundamental group $\pi_1^{t, ab}(X)$ of $X$. By definition, this group classifies finite abelian Galois covers of $\fX$ which are \'etale over $X$ and tamely ramified at codimension 1 points of $\fX\setminus X$ (i.e. the ramification index at such a point is prime to $p$ and the extension of its residue field is separable). One has a direct sum decomposition
$$
\pi_1^{t,ab}(X)\cong\pi_1^{ab}(X)(p')\oplus\pi_1^{ab}(\fX)(p)
$$
where the symbols $(p')$ and $(p)$ stand for the maximal profinite prime-to-$p$ (resp. $p$) quotients of the groups in question (for the $p$-part, notice that any abelian cover of $\fX$ which is of $p$-power degree, \'etale over $X$ and tamely ramified in codimension 1 must be \'etale in codimension 1, hence \'etale by Zariski-Nagata purity). Since the fundamental group is a birational invariant of projective varieties, the above decomposition shows that $\pi_1^{t,ab}(X)$ depends only on $X$ but not on the compactification $\fX$. Needless to say, all these notions and facts are valid more generally over any perfect base field in place of $k$. 

\begin{proposition}
\label{proposition:prop1bis} Under the assumptions of the previous proposition, there is an exact sequence 
\begin{equation}
\label{eqn:tame}
0 \lra T\lra \pt( X) \lra T(Alb_{X}) \lra 0
\end{equation}
where $T(Alb_X)$ denotes the full Tate module of $Alb_X$ and $T$ is a finite abelian group whose prime-to-$p$ part is isomorphic to that of the (dual of) the group $S$ considered above and whose $p$-part is isomorphic to the (dual of) the $p$-primary torsion subgroup of $NS(\fX)$.
\end{proposition}

\noindent{\em Proof.} We use the decomposition of $\pt (X)$ recalled above. The assertion for the prime-to-$p$ part follows from the previous proposition by dualising and passing to the limit. For the $p$-part we note that $T_p(Alb_X)\cong T_p(Alb_{\fX})$, the toric part of $Alb_X$ having no $p$-primary torsion, so the result follows from the analogous statement for $\fX$ proven in (\cite{kala}, Lemma 5).
\enddemo

\section{The Generalisation of Roitman's Theorem}

Keeping the assumptions of the previous section, we now prove Theorem \ref{theorem:thm1}.
We begin by some preliminary observations.
For any positive integer $n$ prime to $p$ the long exact sequence 
\[
\ldots \lra h_i(X) \stackrel{n}{\lra} h_i(X) \lra h_i(X,\bZ/n)\lra h_{i-1}(X)\stackrel{n}{\lra} \ldots
\]
yields a surjection 
\begin{equation}\label{eqn:surj}
h_1(X,\bZ/n)\lra {h_0(X)}_n.
\end{equation}
On the other hand, we have a chain of isomorphisms
\begin{equation}\label{eqn:sv}
h_1(X,\bZ/n)\cong \Hom(h^1(X,\bZ/n),\bZ/n)\cong\Hom(H^1_{\et}(X,\bZ/n),\bZ/n),
\end{equation}
the first by the very definition of the groups in question (note that $\bZ/n$ is injective as a $\bZ/n$-module) and the second by the main result of \cite{sv} (which also holds in positive characteristic by de Jong's work on alterations; compare \cite{ge}). 
Now the key result is:

\begin{proposition}\label{prop:diag}
For any positive integer $n$ prime to $p$ we have a commutative diagram 
\begin{equation}\label{diag} 
\begin{CD}
h_1(X, \bZ/n) @>>> {h_0(X)}_n\\
@V\cong VV @VV alb_X V\\
\Hom(H^1_{\et}(X, \bZ/n), \bZ/n) @>>> {Alb_X(k)}_n
\end{CD}
\end{equation}
where the bottom horizontal map comes from pulling back covers of $Alb_X$ to $X$ and the other unnamed maps come from (\ref{eqn:surj}) and (\ref{eqn:sv}).
\end{proposition}

The proposition yields an immediate proof of Theorem \ref{theorem:thm1} as follows. It is enough to consider $\ell$-primary torsion for a prime $\ell\neq p$. By making $n$ vary among powers of $\ell$ and passing to the direct limit, we get a diagram in which the upper horizontal map is surjective, while those on the left and at the bottom are isomorphisms (the latter by Corollary \ref{corollary:cor1}). Hence all maps in the diagram become isomorphisms in the limit.\medskip

For the proof of the proposition we need the following technical statements about abelian groups whose formal proof will be left to the readers.

\begin{lemma}\label{form}${}{}{}$
\begin{enumerate}
\item For any abelian group $A$ and integer $n>0$ there is a canonical isomorphism $$\Hom(A_n,\bZ/n)\cong \Ext^1(A,\bZ/n).$$
\item Let $(C_{\bullet}, d)$ be a homological complex of free abelian groups concentrated in nonnegative degrees. Then the natural map 
$$
H_1(C_{\bullet}\otimes\bZ/n)\to H_0(C_{\bullet})_n
$$
coming from tensoring by the exact sequence $0\to\bZ\to\bZ\to\bZ/n\to\nobreak 0$ can be identified with the map
\begin{equation}\label{tor}
H_1(C_{\bullet}\otimes\bZ/n)\to \Tor(H_0(C_{\bullet}),\bZ/n)
\end{equation}
coming from computing the Tor-group by means of the free resolution $d(C_1)\to C_0$ of $H_0(C_{\bullet})$. 
\item With the previous notations, the natural map
$$
\Ext^1(H_0(C_{\bullet}),\bZ/n)\to\Ext^1(C_{\bullet},\bZ/n)
$$
induced by the truncation map $C_{\bullet}\to H_0(C_{\bullet})$
can be identified (using statement 1. and the self-injectivity of the ring $\bZ/n$) with the image of the map (\ref{tor}) under the functor $\Hom(\,\, ,\bZ/n)$.
\end{enumerate}
\end{lemma}

\noindent{\it Proof of the Proposition.} We prove the commutativity of the dual diagram
\[
\begin{CD}
\Hom(Alb_X(k)_n, \bZ/n) @>>>H^1_{\et}(X, \bZ/n)\\
@VVV @VV\cong V\\
\Hom({h_0(X)}_n,\bZ/n) @>>> h^1(X, \bZ/n)
\end{CD}
\]
which, using Lemma \ref{form} (1), can be rewritten as
\[
\begin{CD}
\Ext^1(\widetilde{Alb}_X(k), \bZ/n) @>>>H^1_{\et}(X, \bZ/n)\\
@VVV @VV\cong V\\
\Ext^1({h_0(X)},\bZ/n) @>>> \Ext^1(C_{\bullet}(X), \bZ/n)
\end{CD}
\]
where the Ext-groups are taken with respect to the category of abelian groups (there was no harm in replacing $Alb_X$ by $\widetilde{Alb}_X$ since $\bZ$ is torsion-free). Using Lemma \ref{form} the bottom horizontal map can then be identified as coming from the natural truncation map. Now we apply the rigidity theorem of Suslin-Voevodsky (\cite{sv}, Theorem 4.5) to the three Ext-groups and a standard comparison theorem to the fourth group to obtain a diagram  
\begin{equation}\label{diag2}
\begin{CD}
\Ext^1_{\et}(\widetilde{Alb}_X, \bZ/n) @>>>\Ext^1_{\et}(\bZ(X), \bZ/n)\\
@VVV @VV\cong V\\
\Ext^1_{\et}(H_0(C_{\bullet}(\bZt(X)))_{\et},\bZ/n) @>>> \Ext^1_{\et}(C_{\bullet}(\bZt(X)), \bZ/n)
\end{CD}
\end{equation}
where the Ext-groups are now taken on the \'etale site of $Sm/k$, the subscript {\em \'et} means sheafification for the \'etale topology and $\bZ(X)$ is the \'etale sheaf whose sections over a smooth $k$-scheme $Y$ are given by the free abelian group with basis $\Hom(Y,X)$. Note that the rigidity theorem was applicable to the upper left group by virtue of Lemma \ref{transfers} and to the two lower ones by (\cite{sv}, Corollary 7.5). Now to finish the proof, we claim that the above diagram is induced by applying the functor $\Ext^1_{\et}(\,\,\,,\bZ/n)$ to the commutative diagram of complexes of sheaves 
\[
\begin{CD}
C_{\bullet}(\bZt(X)) @>>>H_0(C_{\bullet}(\bZt(X)))_{\et}\\
@AAA @VVV\\
\bZt(X) @>>>\widetilde{Alb}_X
\end{CD}
\]
whose existence was established in Section 3 (the map on the left inducing the inverse of the isomorphism marked in (\ref{diag2})) . 

The identification of the bottom horizontal and left vertical arrows in (\ref{diag2}) follows from the functoriality of the rigidity isomorphism. As for the upper horizontal map, note first that it is well known to be induced by the map of \'etale sheaves $\bZ(X)\to\widetilde{Alb}_X$ which factors through the natural inclusion $\bZ(X)\to\bZt(X)$ by Lemma \ref{transfers}. Now by (\cite{sv}, Corollary 10.10) the natural map $\Ext^1_{\et}(\bZt(X),\bZ/n)\to\Ext^1_{\et}(\bZ(X),\bZ/n)$ can be identified with the map $\Ext^1_{qfh}(\bZt(X)_{qfh},\bZ/n)\to \Ext^1_{qfh}(\bZ(X)_{qfh},\bZ/n)$, where the subscript $qfh$ denotes sheafification for the so-called $qfh$-topology introduced in {\em loc. cit.}, which is finer than the \'etale topology. But the latter map is an isomorphism, for by (\cite{sv}, Theorem 6.7) $\bZt(X)$ can be identified, after localisation by the characteristic $p$, with $\bZ(X)_{qfh}$. This finishes the identification of the upper horizontal map, and for the right vertical map one first uses the same argument to pass from $\bZ(X)$ to $\bZt(X)$, whereupon the result follows from the construction of the isomorphism $h^1(X, \bZ/n)\to H^1_{\et}(X, \bZ/n)$ in the proof of Theorem 7.5 in \cite{sv} (from which one sees that it can be identified with the map induced by $\bZt(X) \to C_*(\bZt(X))$ on $\Ext^1$-groups).  
\enddemo

\begin{remark}\rm
One can give a quicker, albeit less conceptual proof of Theorem \ref{theorem:thm1} which avoids the verification of commutativity, in the spirit of the simplified version of Bloch's approach to Roitman's theorem given in \cite{colliot}. Indeed, using the lower horizontal and left vertical arrows in (\ref{diag}) and passing to the limit one gets a surjection ${Alb_X(k)}_{\ell-tor}\to {h_0(X)}_{\ell-tor}$. Since $Alb_X$ is a semiabelian variety, both groups here must be isomorphic  to some finite direct power of ${\bQ}_{\ell}/{\bZ}_{\ell}$, so that for any $m>0$ the groups ${Alb_X(k)}_{\ell^m}$ and ${h_0(X)}_{\ell^m}$ are finite and the order of the second one doesn't exceed that of the first. On the other hand, one proves by cutting $\fX$ successively by well-chosen hyperplanes that there is a smooth closed curve $C\subset X$ such that the map $Alb_C(k)_{\ell^m}\to Alb_X(k)_{\ell^m}$ is surjective for any $m$ (this is entirely classical in the projective case and the general case follows from it by a localisation argument), whence the surjectivity of $alb_X$ on $\ell^m$-torsion by reduction to the case of curves treated in \cite{lich} and (\cite{sv}, Theorem 3.1). Using the relation between the orders of the $\ell^m$-torsion subgroups this implies that $alb_X$ is an isomorphism when restricted to prime-to-$p$ torsion subgroups.
\end{remark}

\section{Proof of Theorem \ref{theorem:thm2}}

Assume now that $k$ is the algebraic closure of a finite field $\bF$ and denote by $G$ the Galois group $Gal(k|\bF)$. By extending $\bF$ if necessary we may assume that there are varieties $X_{\bF}\subset \fX_{\bF}$ defined over $\bF$ such that $X_{\bF}$ has an $\bF$-rational point and $X_{\bF}\times_{\bF} k\cong X$, $\fX_{\bF}\times_{\bf} k\cong \fX$. Similarly to section 10 of \cite{kasa}, the key to the proof of Theorem \ref{theorem:thm2} is the exact sequence (\ref{eqn:tame}) which in this case is in fact an exact sequence of $G$-modules.  

Recall from Section 4 that the abelianised tame fundamental group can be defined for $X_{\bF}$ as well. Moreover, there is a natural projection $\pt(X_{\bF}) \to G\cong \widehat{\bZ}$ whose kernel $\pt(X_{\bF})^0$ can be identified with the coinvariants of $\pt(X)$ under the action of $G$. Therefore taking coinvariants under Frobenius in the exact sequence (\ref{eqn:tame})
yields the sequence
\[
0\to T_G\lra \pt(X_{\bF})^0\lra  Alb_{X_{\bF}}(\bF) \lra 0.
\]
By the main result of \cite{ss} the group in the middle is isomorphic to ${h_0(X_{\bF})}^0$ by means of a reciprocity map $h_0(X_{\bF})^0\to\pt(X_{\bF})$ which sends the class of a closed point of $X$ to the class of its Frobenius. Using this isomorphism and taking the limit of the direct system of similar exact sequences over finite extensions $\bF'$ of $\bF$ we get an isomorphism 
\[
h_0(X)^0\cong Alb_{X}(k)
\]
because the direct limit of the finite groups $T_{Gal(k|{\bF'})}$ is trivial. It remains to see that this isomorphism is induced by the Albanese map. For this it suffices to consider the image of the class of a zero-cycle of the form $P_1-P_2$ and we may replace $k$ by the finite extension of $\bF$ over which both $P_1$ and $P_2$ are defined. Moreover, by using the covariant functoriality of the Albanese map and of the reciprocity map we may assume that $P_1$ and $P_2$ both lie on some smooth curve $C$ and check the required compatibility for $C$, but this is a well-known property of Lang's class field theory (see \cite{serre2}).

Finally we note that the above proof has the following interesting by-product, generalising the similar statement proved in \cite{kasa} for the proper case:

\begin{corollary}
With notations as above, the natural map $h_0(X_{\bF})\to h_0(X)$ has a finite kernel isomorphic to the group $T$ introduced in Proposition \ref{proposition:prop1bis}.
\end{corollary}

\bigskip
\begin{tabbing}
 \hspace{1cm}\= Michael Spie{\ss}\hspace{5cm}\= Tam\a'as Szamuely \\
 \>School of Mathematical Sciences \> Alfr\a'ed R\a'enyi Institute of Mathematics\\
 \>University of Nottingham \> Hungarian Academy of Sciences\\
  \>University Park \> PO Box 127 \\ 
 \>NG7 2RD Nottingham \> H-1364 Budapest   \\ 
 \>United Kingdom \> Hungary     
 \end{tabbing}

\noindent \hspace{.9cm} e-mail: mks@maths.nott.ac.uk \\ 
\hspace*{2.3cm}  szamuely@renyi.hu

\end{document}